# EXAMPLES OF NON-TRIVIAL CONTACT MAPPING CLASSES IN ALL DIMENSIONS

PATRICK MASSOT AND KLAUS NIEDERKRÜGER

ABSTRACT. We give examples of contactomorphisms in every dimension that are smoothly isotopic to the identity but that are not contact isotopic to the identity. In fact, we prove the stronger statement that they are not even symplectically pseudo-isotopic to the identity. We also give examples of pairs of contactomorphisms which are smoothly conjugate to each other but not by contactomorphisms.

## 1. INTRODUCTION

In our paper [MNW13] with Chris Wendl, we described examples of high dimensional contact manifolds that behave in many ways similarly to the tight contact structures on the 3-torus. As explained in [GM15], the tight contact structures on the 3-torus with positive Giroux torsion admit contactomorphisms that are smoothly isotopic to the identity but not through contactomorphisms[1]. The goal of this note is to observe that our high dimensional examples also admit such contactomorphisms. But let us first recall briefly the situation in dimension three. On $\mathbb{T}^3 = (\mathbb{R}/2\pi\mathbb{Z})^3$ with coordinates $(s, t, \theta)$, consider the family of contact structures $\xi_n$ with $n \in \mathbb{Z}_{>0}$ defined by $\xi_n := \ker(\sin(ns)\,dt + \cos(ns)\,d\theta)$. The rotation maps

$$\Psi_{n,m} \colon (\mathbb{T}^3, \xi_n) \to (\mathbb{T}^3, \xi_n),\, (s,t,\theta) \mapsto \left(s + \tfrac{2\pi m}{n}, t, \theta\right),$$

(with $0 \leq m < n$) are contactomorphisms which are smoothly isotopic to the identity but, according to [GM15], they are pairwise non-isotopic through contactomorphisms.

We now turn to higher-dimensional examples. Let $\mathbf{k} \subset \mathbb{R}$ be a field of real numbers such that $\dim_{\mathbb{Q}} \mathbf{k}$ is finite and such that $\mathbf{k}$ is totally real (i.e. any field embedding $\mathbf{k} \hookrightarrow \mathbb{C}$ is real-valued). In [MNW13, Theorem 9.10] we associated to $\mathbf{k}$ a compact manifold $M_{\mathbf{k}}$ equipped with 1-forms $\alpha_\pm$ such that the formula

$$\xi_n := \ker\left(\frac{1 + \cos(ns)}{2}\alpha_+ + \frac{1 - \cos(ns)}{2}\alpha_- + \sin(ns)\,dt\right)$$

for $n \geq 1$ defines a family of contact structures on $\mathbb{T}^2 \times M_{\mathbf{k}}$, where $(s, t)$ are the coordinates of $\mathbb{T}^2$. The proof of Observation 4.1 contains a short explanation on how the $M_{\mathbf{k}}$ arise. The contact manifolds $(\mathbb{T}^2 \times M_{\mathbf{k}}, \xi_n)$ have the following properties:
- They all admit Reeb vector fields without contractible closed orbits.
- They are all homotopic through almost contact structures but not contactomorphic.
- $\xi_n$ is strongly fillable only if $n = 1$.

For instance $M_{\mathbb{Q}} = \mathbb{S}^1$ with $\alpha_\pm = \pm d\theta$ so $\mathbb{T}^2 \times M_{\mathbb{Q}} = \mathbb{T}^3$ with $\xi_n$ as above.

Note that there are infinitely many such fields $\mathbf{k}$ for each given $\dim_{\mathbb{Q}} \mathbf{k} > 1$, and the corresponding $M_{\mathbf{k}}$ are pairwise non-homeomorphic.

**Theorem 1.1.** *For any totally real number field $\mathbf{k}$, any $n$ greater than one and any $1 \leq m < n$, the contactomorphism*

$$\Psi_{n,m} \colon (\mathbb{T}^2 \times M_{\mathbf{k}}, \xi_n) \to (\mathbb{T}^2 \times M_{\mathbf{k}}, \xi_n),\, (s,t,\theta) \mapsto \left(s + \tfrac{2\pi m}{n}, t, \theta\right)$$

*is smoothly isotopic to the identity but it is not symplectically pseudo-isotopic to the identity, so in particular it is not contact isotopic to the identity. In addition, there is a contactomorphism which is conjugated to $\Psi_{n,m}$ inside $\mathrm{Diff}(\mathbb{T}^2 \times M_{\mathbf{k}})$ but not inside $\mathrm{Diff}(\mathbb{T}^2 \times M_{\mathbf{k}}, \xi_n)$.*

---
[1]To our knowledge such contactomorphisms were first exhibited by Gompf on $\mathbb{S}^1 \times \mathbb{S}^2$, see [Gom98].

The 3-dimensional result mentioned above has been obtained in [GM15] using Giroux's theory of $\xi$-convex surfaces. Such methods do not seem to be sufficiently powerful to prove the higher dimensional results treated in this text, and even in dimension 3 it seems unlikely that they might yield the stronger pseudo-isotopy obstruction. Instead, we will use $J$-holomorphic curve techniques to show that a certain pre-Lagrangian submanifold $P$ in $\mathbb{T}^2 \times M_\mathbf{k}$ cannot be displaced from itself by any contactomorphism that is symplectically pseudo-isotopic to the identity. The main theorem follows because $\Psi_{m,n}$ does displace $P$.

**Outline.** In Section 2 we explain how soft methods reduce the non-displaceability statement to a statement about non-existence of weakly exact closed Lagrangian submanifolds in $S\xi \times \mathbb{C}$ where $S\xi$ is the symplectization of a closed hypertight contact manifold. If $S\xi$ were geometrically bounded, this would just be a special case of a result by Gromov [Gro85, Section 2.3.$B_3'$], here we need to combine it with Hofer's compactness for holomorphic disks in symplectizations [Hof93]. This requires some care because the end of $S\xi \times \mathbb{C}$ is neither convex nor concave, and because neither the closed Lagrangian submanifold serving as boundary condition for an inhomogeneous Cauchy-Riemann problem, nor the perturbation term involved are in product form. We explain the solution to these problems in Section 3 in detail, even if the resulting proofs are essentially classical. In Section 4 we apply the non-displaceability result to our examples to obtain Theorem 1.1.

**About technological sophistication.** Note that a stronger non-displaceability result holds: The pre-Lagrangian $P$ contains a Legendrian submanifold $\Lambda$ which cannot be disjoined from $P$. This can be proved by setting up a Floer theory for Lagrangian lifts of $P$ and $\Lambda$ in the symplectization of $\xi_n$ as was done in [EHS95] (see Lemma 2.4 about why invariance under compactly supported Hamiltonian isotopies is enough). Such a strategy involves a lot more technical work than is necessary to deduce our theorem on contact transformations. An even more high-tech road would be to prove that contact transformations which are symplectically pseudo-isotopic to the identity act trivially on contact homology and use it to prove Theorem 1.1. However we feel that such a monumental proof would not make sense as long as our only examples can be handled by much more elementary techniques. So we chose instead to prove the weaker non-displaceability result (which is also of independent interest and has less hypotheses). Here one can also envision variations on the argument. One referee pointed out to us that we could adapt to our setup the variation on Gromov's argument which is explained in [MS04, end of Section 9.2]. This variation uses holomorphic strips instead of disks and is arguably slightly more contrived but does not set up a full Floer theory so it is also elementary in the sense of the current discussion. Note however that such a road would bypass Theorem 2.6 which has independent interest.

**Acknowledgments.** The authors would like to thank Kyler Siegel for asking us about examples of contactomorphisms that are conjugated inside the diffeomorphism group but not inside the contactomorphism group. We are also indebted to Chris Wendl and Baptiste Chantraine for valuable discussions and to the anonymous referees for suggesting improvements to the exposition. We are especially grateful to Sylvain Courte who pointed out to us that our original proof of Observation 2.1 was not correct. Klaus Niederkrüger is thankful for the kind invitation to the Centre de Mathématiques Laurent Schwartz of the École polytechnique during which this project was started. Klaus Niederkrüger is currently being supported by the ERC Advanced Grant LDTBud.

## 2. From pseudo-isotopic disjunctions to weakly exact Lagrangians

In this section, we discuss symplectic pseudo-isotopies, and pre-Lagrangian submanifolds. We then explain how persistence of certain Lagrangian intersections in the symplectization implies persistence of pre-Lagrangian intersections in the corresponding contact manifold. Finally we explain how existence of relevant Lagrangian intersections follows from a result about weakly exact Lagrangians which will be proved in the next section.

First recall that the symplectization of a contact manifold $(M, \xi)$ is

$$S\xi = \{p \in T^*M \mid \ker p = \xi\} \xrightarrow{\pi} M \ ,$$

equipped with its canonical Liouville form $\lambda$ and symplectic structure $\omega = d\lambda$ and its $\mathbb{R}$-action $\tau_s(p) := e^s p$. In this paper all contact structure are cooriented and $\ker p = \xi$ is meant as an equality of cooriented hyperplanes so that $S\xi$ is diffeomorphic to $\mathbb{R} \times M$ with its obvious $\mathbb{R}$-action. We denote by $S_{-\infty}\xi$ its negative end and by $S_{+\infty}\xi$ its positive end. By definition, a neighborhood of $S_{\pm\infty}\xi$ is a set containing an open set which is invariant under the action of $\mathbb{R}_\pm$.

**Definition** (Cieliebak-Eliashberg [CE12, Section 14.5]). A **symplectic pseudo-isotopy** $F$ of a contact manifold $(M, \xi)$ is a symplectomorphism
$$F \colon (S\xi, \omega) \to (S\xi, \omega)$$
that restricts on a neighborhood of the negative end $S_{-\infty}\xi$ to the identity, and that preserves the Liouville form $\lambda$ on a neighborhood of the positive end.

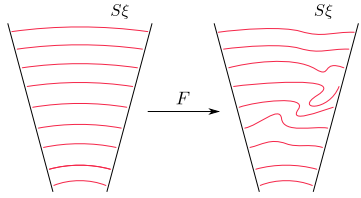

FIGURE 1. A pseudo-isotopy commutes with the $\mathbb{R}$-action on neighborhoods of $S_{\pm\infty}\xi$, but not necessarily in between.

Because a symplectic pseudo-isotopy $F$ preserves the canonical 1-form $\lambda$ on some $\mathbb{R}_+$-invariant neighborhood $U_+$ of $S_{+\infty}\xi$, it also preserves the vector field $X$ generating the $\mathbb{R}$-action that is characterized by $d\lambda(X, \cdot) = \lambda$. Hence it commutes on $U_+$ with the $\mathbb{R}_+$-action and induces a contactomorphism of $(M, \xi)$ as follows: for any $x$ in $M$ choose any $p$ in $U_+$ above $x$ and define $\varphi(x) = \pi(F(p))$. More concretely, any contact form $\alpha$ for $\xi$ is a section of $S\xi$ which identifies $S\xi$ with:
$$\left(\mathbb{R} \times M, \omega = d(e^t \alpha)\right),$$
where $t$ denotes the coordinate on the $\mathbb{R}$-factor. Commutation with the $\mathbb{R}_+$-action on $U_+$ means there is a function $g \colon M \to \mathbb{R}$ and a diffeomorphism $\varphi \colon M \to M$ such that $F(t, x) = (t - g(x), \varphi(x))$ for $t \gg 0$. Since $F$ preserves $\lambda = e^t \alpha$ in this region, we get that $\varphi^* \alpha = e^g \alpha$.

A contactomorphism obtained this way is said to be **symplectically pseudo-isotopic to the identity**. These contactomorphisms form a group denoted by[2] $\mathrm{Diff}_\mathcal{P}(M, \xi)$.

**Observation 2.1** (Implicit in [CE12]). *The group $\mathrm{Diff}_\mathcal{P}(M, \xi)$ contains the identity component of $\mathrm{Diff}(M, \xi)$, ie the group of contactomorphisms which are contact isotopic to the identity.*

*Proof.* Let $\varphi_s \colon M \to M$ with $s \in [0, 1]$ be a contact isotopy and let $X_s$ be the vector field generating $\varphi_s$ (characterized by $d\varphi_s/ds = X_s \circ \varphi_s$). The Hamiltonian function $H_s \colon S\xi \to \mathbb{R}$ sending $p$ to $p(X_s)$ defines a Hamiltonian isotopy $\hat{\varphi}_s$ of $S\xi$ which is an $\mathbb{R}$-equivariant lift of $\varphi_s$. We will construct a symplectic pseudo-isotopy inducing $\varphi_1$ by cutting off $\hat{\varphi}_1$. Let $U_+$ be an $\mathbb{R}_+$-invariant closed neighborhood of $S_{+\infty}\xi$ which is disjoint from some neighborhood of $S_{-\infty}\xi$. Let $U_{++}$ be a neighborhood of $\cup_{s \in [0,1]} \hat{\varphi}_s(U_+)$ and let $\rho$ be a smooth function from $S\xi$ to $[0, 1]$ which vanishes on a neighborhood of $S_{-\infty}\xi$ and equals 1 on $U_{++}$. Let $\psi_s$ be the Hamiltonian flow generated by $\rho H_s$. This flow exists for all $s$ in $[0, 1]$ because $\rho H_s$ is equivariant outside a compact subset of $S\xi$. And it coincides with $\hat{\varphi}_s$ on $U_+$ for all $s$ by construction of $U_{++}$. Its time one map $\psi_1$ is the desired pseudo-isotopy. □

We now turn our attention to pre-Lagrangian submanifolds. For any submanifold $P \xhookrightarrow{\iota} (M, \xi)$, a **lift** of $P$ is a submanifold $L$ in $S\xi$ which is transverse to the $\mathbb{R}$-action and projects onto $P$. Any lift of $P$ can be seen as $\alpha(P)$ for some contact form $\alpha$. Because the Liouville form $\lambda$ has the tautological property $\alpha^* \lambda = \alpha$ for any $\alpha$, we get that $\alpha(P)$ is an isotropic submanifold of $S\xi$ if and only if $\iota^* \alpha$ is closed. This motivates the following definition (attributed to Bennequin in [EHS95]).

---

[2]This is not exactly the definition written in [CE12] but it is what the authors intended to write.

**Definition.** Let $(M, \xi)$ be a $(2n + 1)$-dimensional contact manifold. An $(n + 1)$-dimensional submanifold $P \subset M$ is **pre-Lagrangian** if there is a contact form for $\xi$ whose restriction to $P$ is closed.

Note that a pre-Lagrangian submanifold $P$ is always transverse to $\xi$ since otherwise its tangent space at a non-transverse point would be an $(n+1)$-dimensional isotropic subspace of $\xi$.

Remember that a Lagrangian $L$ in a symplectic manifold $(W, \omega)$ is called **weakly exact** if $\int_{\mathbb{D}^2} u^*\omega$ vanishes for every smooth map $u \colon (\mathbb{D}^2, \partial\mathbb{D}^2) \to (W, L)$. On the pre-Lagrangian side, one can prove the following:

**Lemma 2.2.** *Let $P \xhookrightarrow{\iota} (M, \xi)$ be a pre-Lagrangian submanifold and denote by $\mathcal{A}(P)$ the space of contact forms for $\xi$ whose restrictions to $P$ are closed. If $P$ is closed then the following properties are equivalent:*
  (1) *there exists an $\alpha$ in $\mathcal{A}(P)$ such that $\int_{\mathbb{D}^2} u^* d\alpha$ vanishes for every smooth map $u \colon (\mathbb{D}^2, \partial\mathbb{D}^2) \to (M, P)$.*
  (2) *for every $\alpha$ in $\mathcal{A}(P)$, $\int_{\mathbb{D}^2} u^* d\alpha$ vanishes for every smooth map $u \colon (\mathbb{D}^2, \partial\mathbb{D}^2) \to (M, P)$.*
  (3) *there is a Lagrangian lift of $P$ which is weakly exact*
  (4) *all Lagrangian lifts of $P$ are weakly exact.*

A closed pre-Lagrangian submanifold with any of the above properties will be called **weakly exact**. This terminology parallels the Lagrangian case but note that there is nothing like a strongly exact closed pre-Lagrangian (see the proof below).

*Sketch of proof of Lemma 2.2.* For any $\alpha$ in $\mathcal{A}(P)$, the restriction $\iota^*\alpha$ cannot be exact since $P$ is closed so it does not have a nowhere vanishing exact 1-form. Also Tischler's theorem proves that the leaves of the foliation printed by $\xi$ on $P$ are either dense or coincide with the fibers of a fibration $P \to \mathbb{S}^1$. Analyzing both cases, one can show that cohomology classes of $\iota^*\alpha$ for $\alpha$ in $\mathcal{A}(P)$ are all positively proportional (see [EHS95, Proposition 2.2.2] for details about that argument). The rest is then an easy consequence of Stokes' theorem, the fact that both the projection $\pi \colon S\xi \to M$ and the sections $\alpha$ are homotopy equivalences and the tautological property of $\lambda$. $\square$

In the rest of this section we will explain how the following result, which is the key step to the proof of our main theorem, can be translated into a statement about the non-existence of certain weakly exact Lagrangians. Remember that a contact structure is called **hypertight** if it admits a Reeb vector field without contractible closed orbits.

**Theorem 2.3.** *A closed weakly exact pre-Lagrangian in a closed hypertight contact manifold cannot be displaced by any contactomorphism that is symplectically pseudo-isotopic to the identity.*

Note that one indeed needs the weak exactness assumption in the theorem since Darboux balls are displaceable by contact isotopy and contain plenty of closed pre-Lagrangian submanifolds.

**Lemma 2.4.** *Let $\varphi$ be a contactomorphism of $(M, \xi)$ that is symplectically pseudo-isotopic to the identity. For any compact subset $K$ in $S\xi$, there is a compactly supported Hamiltonian isotopy $\Phi$ in $S\xi$ such that, for every $p \in K$,*
$$\pi(\Phi_1(p)) = \varphi(\pi(p)),$$
*where $\pi \colon S\xi \to M$ is the canonical projection.*

*Proof.* Let $F$ be a symplectic pseudo-isotopy inducing $\varphi$. By definition of symplectic pseudo-isotopies, there are $\mathbb{R}_\pm$-invariant neighborhoods $U_-$ of $S_{-\infty}\xi$ and $U_+$ of $S_{+\infty}\xi$ respectively such that
  - $F|_{U_-} = \mathrm{Id}_{U_-}$;
  - $F|_{U_+}$ commutes with the action $\tau_s$ for any positive $s$.

For any $s$ in $\mathbb{R}$ the conjugation $\tau_s \circ F \circ \tau_{-s}$ is also a symplectic pseudo-isotopy inducing $\varphi$ and has support outside $\tau_s(U_-)$. We replace $F$ by $\tau_s \circ F \circ \tau_{-s}$ for a sufficiently large positive $s$, to ensure that the compact subset $K$ is entirely contained in $U_-$.

We next consider the symplectic isotopy $\Phi$ given by the commutators
$$\Phi_s := [\tau_{-s}, F].$$
The support of each $\Phi_s$ lies outside $F(U_+) \cup \tau_{-s}(U_-)$ hence is compact. It follows that $\Phi$ is a Hamiltonian isotopy. Indeed, if $X_s$ is a time-dependent vector field generating $\Phi$, each 1-form $\iota_{X_s}\omega$ is closed and has compact support. Since the inclusion of any level $\{t\} \times M$ into $\mathbb{R} \times M$ is a homotopy equivalence, it induces an isomorphism in de Rham cohomology so that compactly supported closed forms are all exact.

Let $p$ be any point in $K$. Since $K$ is in $U_-$, $\Phi_s(p) = \tau_{-s} \circ F \circ \tau_s(p)$ for all $s$. Let $s_1$ be a positive number so large that $\tau_{s_1}(K)$ is contained in $U_+$. Then
$$\pi(\Phi_{s_1}(p)) = \pi(F(\tau_{s_1}(p))) = \varphi(\pi(p)).$$
Of course we can replace $s \mapsto \Phi_s$ by $s \mapsto \Phi_{s/s_1}$ to make sure that $s_1 = 1$. $\square$

**Proposition 2.5** ([Gro85, Section 2.3.$B'_3$]). *Let $(W, d\lambda)$ be an exact symplectic manifold, and let $L \subset W$ be a Lagrangian. If $\Phi_s \colon W \to W$ is a Hamiltonian isotopy, then we can find a Lagrangian immersion*
$$j \colon L \times \mathbb{S}^1 \looparrowright (W \times \mathbb{C}, d\lambda \oplus dx \wedge dy)$$
*where $x + iy$ is the coordinate on $\mathbb{C}$, such that the self-intersection points of $j(L \times \mathbb{S}^1)$ are in one-to-one correspondence with the intersection points in $L \cap \Phi_1(L)$. If $j$ is an embedding and if $L$ is weakly exact, then $j$ will also be weakly exact.*

A detailed proof of this proposition can be found in [ALP94, Theorem 2.3.6]. In the next section we will combine ideas by Gromov [Gro85, Section 2.3.$B'_3$] with the compactness in [Hof93] by Hofer to prove the following theorem.

**Theorem 2.6.** *If $(M, \xi)$ is a closed contact manifold that is hypertight, then $(S\xi \times \mathbb{C}, d\lambda \oplus dx \wedge dy)$ does not contain any weakly exact closed Lagrangian.*

Using all this we can prove Theorem 2.3. Suppose that $P$ is a closed weakly exact pre-Lagrangian submanifold in a hypertight $(M, \xi)$. Let $\varphi$ be a contactomorphism symplectically isotopic to the identity and let $L_P$ be a Lagrangian lift of $P$. According to Lemma 2.2, $L_P$ is weakly exact. Assume for contradiction that $P \cap \varphi(P) = \varnothing$. Lemma 2.4 applied to $K = L_P$ and $\varphi$ gives a Hamiltonian isotopy $\Phi$ in $(S\xi, d\lambda)$ which displaces $L_P$: $L_P \cap \Phi_1(L_P) = \varnothing$. Proposition 2.5 turns it into a weakly exact embedded Lagrangian in $(S\xi \times \mathbb{C}, d\lambda \oplus dx \wedge dy)$, which contradicts Theorem 2.6.

3. From hypertightness to absence of weakly exact Lagrangians

In this section we prove Theorem 2.6. following Gromov's argument in [Gro85, Section 2.3.$B'_3$]. The strategy is to show that there is a non-trivial holomorphic disk with boundary on any closed Lagrangian submanifold of $S\xi \times \mathbb{C}$. These disks result from bubbling of an inhomogeneous Cauchy-Riemann equation.

We fix a contact form $\alpha$ without contractible Reeb orbit. We identify $(S\xi, d\lambda)$ with $(\mathbb{R} \times M, d(e^t\alpha))$ using the contact form $\alpha$ and denote by $\pi_\xi$, $\pi_\mathbb{R}$, $\pi_M$ and $\pi_\mathbb{C}$ the canonical projections of $S\xi \times \mathbb{C}$ to $S\xi$, $\mathbb{R}$, $M$ and $\mathbb{C}$ respectively. We fix an $\mathbb{R}$-invariant almost complex structure $J_\alpha$ on $\mathbb{R} \times M$ which preserves $\xi$, is compatible with the restriction of $d\alpha$ to $\xi$ and sends $\partial_t$ to $R_\alpha$. Let $L \subset S\xi \times \mathbb{C}$ be a closed Lagrangian, $U_L$ a compact tubular neighborhood of $L$ and $p_0$ a point in $L$. We assume that $\pi_\mathbb{R}(U_L)$ lies in $\{t > 1\}$ (this can be arranged by a constant rescaling of $\alpha$). All these objects, including $J_\alpha$, are now fixed forever. We denote by $\mathcal{B}$ the space of $W^{1,p}$-maps $u$ from $(\mathbb{D}^2, \partial\mathbb{D}^2, 1)$ to $(S\xi \times \mathbb{C}, L, p_0)$ which are homotopic to the constant map $u_0 \colon z \mapsto p_0$.

We will consider inhomogeneous Cauchy-Riemann equations
$$\bar\partial_J u = G(u),$$
where $J = J_\alpha \oplus i$ on $S\xi \times \mathbb{C}$ and $u \in \mathcal{B}$ is the unknown. The perturbation term $G$ is a section of the following bundle of complex-antilinear maps:
$$\overline{\mathrm{Hom}}_\mathbb{C}(T\mathbb{D}^2, T(S\xi \times \mathbb{C})) \to \mathbb{D}^2 \times (S\xi \times \mathbb{C})$$

and $G(u)\colon \mathbb{D}^2 \to \overline{\mathrm{Hom}}_{\mathbb{C}}(T\mathbb{D}^2, T(S\xi \times \mathbb{C}))$ denotes the restriction of $G$ to the graph of $u$: $G(u)(z) = G(z, u(z))$.

Let $G$ be a family of perturbation terms $G_s$ for $s \in [0, 1]$ and set

(3.1) $$\mathcal{M}(G) = \left\{(s, u) \in [0, 1] \times \mathcal{B} \,\Big|\, \bar{\partial}_J u = G_s(u)\right\}.$$

The spaces of perturbation terms we use are:
$$\mathcal{G}_{\varepsilon, C} = \left\{(\mathbf{0} \oplus Cs\, d\bar{z}) + H_s \,\Big|\, \mathrm{supp}\, H \subset (\varepsilon, 1-\varepsilon) \times \mathbb{D}^2 \times U_L\right\}$$

where $\mathbf{0}$ is the 0-section in $\overline{\mathrm{Hom}}_{\mathbb{C}}(T\mathbb{D}^2, T(S\xi))$ and $Cs$ is in $\mathbb{R} \subset T\mathbb{C}$. The term $H$ is a $C^1$ section of the bundle
$$\overline{\mathrm{Hom}}_{\mathbb{C}}(T\mathbb{D}^2, T(S\xi \times \mathbb{C})) \to [0, 1] \times \mathbb{D}^2 \times (S\xi \times \mathbb{C}).$$

**Proposition 3.1.** *If $\varepsilon > 0$ is chosen sufficiently small and $C$ sufficiently large, then there is some $G$ in $\mathcal{G}_{\varepsilon, C}$ such that $\mathcal{M}(G)$ is a smooth $1$-dimensional manifold whose boundary is $\{(0, u_0)\}$ where $u_0$ is the constant disk at $p_0$.*

*Proof.* Let $G$ be any perturbation term in any $\mathcal{G}_{\varepsilon, C}$. Suppose for contradiction that $u$ is a map in $\mathcal{B}$ such that $(1, u)$ is in $\mathcal{M}(G)$. Then $\pi_{\mathbb{C}} \circ u$ is harmonic and the arguments in [Gro85, Section 2.2.$B_3$] (or [ALP94, Section X.5.3]) give a contradiction if $C$ is bigger than $\mathrm{diam}(\pi_{\mathbb{C}}(L))$. Fix one such $C$, so that the operator $\bar{\partial}_J - G_s$ is in particular trivially transverse to the zero section when $s$ is close to 1.

For $s = 0$, the only $(s, u)$ in $\mathcal{M}(G)$ is $(0, u_0)$ because $u$ has to be a homotopically trivial holomorphic disk. The Riemann-Roch formula and our point constraint at $1 \in \mathbb{D}^2$ prove that the index of the linearization $D_{u_0}$ of $\bar{\partial}_J$ at $u_0$ is zero. Since $u_0$ is constant one can explicitly study $D_{u_0}$ and see that its kernel is trivial. So we get that $D_{u_0}$ is surjective.

Hence we can choose $\varepsilon$ small enough to ensure that the operator $\bar{\partial}_J - G_s$ is transverse to the zero section when $s$ is in $[0, \varepsilon)$ and there is no solution when $s$ is in $(1 - \varepsilon, 1]$.

Next we consider the universal moduli space which is the zero set of $\mathcal{D}\colon (s, u, G) \mapsto \bar{\partial}_J u - G_s(u)$. We only need to prove that this operator has surjective linearization at every solution since Sard-Smale applied to $(s, u, G) \mapsto G$ then gives the desired $G$. The linearization of $\mathcal{D}$ at $(s, u, G)$ operates on a triple $(a, \eta, H)$ where $a$ is any real number, $\eta$ is a $W^{1,p}$ section of $u^*T(S\xi \times \mathbb{C})$ which takes values in $u^*TL$ along $\partial \mathbb{D}^2$ and vanishes at 1 and $H$ is a family of perturbation terms (still of class $C^1$) whose support is in $(\varepsilon, 1-\varepsilon) \times \mathbb{D}^2 \times U_L$. The linearization of $\mathcal{D}$ is given by
$$(a, \eta, H) \mapsto D_u \eta - a\dot{G}_s - \partial_p G_s \eta - H_s$$
where $D_u$ is the linearization of $\bar{\partial}_J$ at $u$, $\dot{G}_s = dG_s/ds$ and $\partial_p G_s$ is the derivative of $G_s$ in the direction of $S\xi \times \mathbb{C}$.

If $s$ is larger than $1 - \varepsilon$ then there is nothing to prove since there is no solution. If $s$ is less than $\varepsilon$ then the last two terms disappear and $D_u$ is surjective so the sum is surjective.

For $s$ in $(\varepsilon, 1-\varepsilon)$ we argue by contradiction. If the cokernel is not trivial then there is some $L^q$ section $\alpha$ of $\overline{\mathrm{Hom}}_{\mathbb{C}}(T\mathbb{D}^2, u^*T(S\xi \times \mathbb{C}))$ such that, for every $\eta$ and $H$,
$$\langle \dot{G}_s, \alpha \rangle_{L^2} = 0$$
$$\langle D_u \eta - \partial_p G_s \eta, \alpha \rangle_{L^2} = 0$$
$$\langle H_s, \alpha \rangle_{L^2} = 0.$$

The operator $D_u - \partial_p G_s$ is a real Cauchy-Riemann type operator. Elliptic regularity for its formal adjoint proves that $\alpha$ is smooth. The similarity principle proves that zeros of $\alpha$ are isolated. Since $u$ maps $\partial \mathbb{D}^2$ to $L$, there is an open set $V$ in $\mathbb{D}^2$ such that $u$ maps $V$ to $U_L$ and, after shrinking $V$, $\alpha$ does not vanish on $V$. For any cut-off function $\rho$ with support in $(\varepsilon, 1-\varepsilon) \times V \times U_L$, $H = \rho\alpha$ is an admissible perturbation term which contradicts the last equation above. $\square$

We now choose one $G$ given by Proposition 3.1 and keep it until the end of this section. Since $\mathcal{M}(G)$ is a 1-dimensional manifold with only one boundary point, it cannot be compact. We want to prove that, under the assumptions of hypertightness of $\xi$ and compactness of $M$ and $L$, the only source of non-compactness for $\mathcal{M}(G)$ is bubbling of holomorphic disks so that $L$ is not weakly exact.

We first note some $C^0$-estimates for all solutions in $\mathcal{M}(G)$. For any $(s,(u_\xi, u_\mathbb{C}))$ in $\mathcal{M}(G)$, the component $u_\xi$ is $J_\alpha$-holomorphic outside the preimage of $U_L$. Hence it cannot enter any neighborhood of $S_{+\infty}\xi$ which is disjoint from the projection $\pi_\xi(U_L)$. Similarly, the component $u_\mathbb{C}$ is harmonic outside the preimage of $U_L$ and this implies that the image of $u_\mathbb{C}$ is contained in a fixed compact subset (any disk around $0$ which contains $\pi_\mathbb{C}(U_L)$ is big enough, see for example the proof of [MS04, Lemma 9.2.3]). Those observations are summarized in the following lemma.

**Lemma 3.2.** *There is a neighborhood $U_+$ of $S_{+\infty}\xi$ and a compact set $K_\mathbb{C} \subset \mathbb{C}$ such that, for all $(s,u)$ in $\mathcal{M}(G)$, $u(\mathbb{D}) \subset (S\xi \setminus U_+) \times K_\mathbb{C}$.* □

Next we need some energy bounds. In view of our later use of Hofer's energy, we will introduce the following class of symplectic forms. We consider the space of probe functions

$$\mathcal{F} := \{\psi \colon \mathbb{R} \to \mathbb{R} \mid \psi \text{ is a smooth embedding and } \psi(t) = t \text{ for } t > 1\}$$

and the associated exact symplectic forms $\omega_\psi := d(e^\psi \alpha)$ on $S\xi$.

**Proposition 3.3.** *There is some bound $A$ such that $\left|\int_V u^*(\omega_\psi \oplus \omega_\mathbb{C})\right| \leq A$ for all (measurable) subsets $V \subset \mathbb{D}$, all $(s,u)$ in $\mathcal{M}(G)$ and all $\psi$ in $\mathcal{F}$.*

*Proof.* The first observation, due to Gromov, is that one can turn the inhomogeneous Cauchy-Riemann problem defining $\mathcal{M}(G)$ into an homogeneous one which allows easier energy estimates. To any $u$ in $\mathcal{B}$ we associate its graph

$$\tilde{u} \colon \mathbb{D}^2 \to \mathbb{D}^2 \times (S\xi \times \mathbb{C}), \ z \mapsto (z, u(z))$$

and for any $s$ in $[0,1]$ we consider the almost complex structure $J_s$ on $\mathbb{D}^2 \times (S\xi \times \mathbb{C})$ given by

$$J_s(\dot{z}, \dot{p}) := (i\dot{z}, J\dot{p} + 2G_s \cdot i\dot{z})$$

for every vector $\dot{z} \in T\mathbb{D}^2$ and $\dot{p} \in T(S\xi \times \mathbb{C})$. The pair $(s,u)$ is in $\mathcal{M}(G)$ if and only if $\tilde{u}$ is a $J_s$-holomorphic map.

**Lemma 3.4.** *If $K > 0$ is large enough then $\tilde{\omega}_\psi = (K\omega_\mathbb{D}) \oplus \omega_\psi \oplus \omega_\mathbb{C}$ tames $J_s$ for all $\psi$ and $s$.*

*Proof.* We denote by $|\cdot|_\psi$ the norm associated to $\omega_\psi(\cdot, J_\alpha \cdot)$. For any $v = v_\mathbb{D} + v_\xi + v_\mathbb{C} \in T(\mathbb{D}^2 \times S\xi \times \mathbb{C})$ we have

$$\tilde{\omega}_\psi(v, J_s v) = K|v_\mathbb{D}|^2 + |v_\mathbb{C}|^2 + |v_\xi|_\psi^2 + (\omega_\psi \oplus \omega_\mathbb{C})((v_\xi, v_\mathbb{C}), 2G_s \cdot iv_\mathbb{D}).$$

Outside $U_L$ we have $G_s = 0 \oplus Cs\,d\bar{z}$ so that

$$\tilde{\omega}_\psi(v, J_s v) \geq K|v_\mathbb{D}|^2 + |v_\mathbb{C}|^2 + |v_\xi|_\psi^2 - 2C|v_\mathbb{D}||v_\mathbb{C}|$$

and we get tameness as soon as $K$ is larger than $C^2$. Inside $U_L$, $\psi(t) = t$ and we set $v_X = v_\mathbb{C} + v_\xi$. In particular $|v_X|^2 = \omega_\mathbb{C}(v_\mathbb{C}, iv_\mathbb{C}) + d\lambda(v_\xi, J_\alpha v_\xi)$ is independent of $\psi$. Since

$$\tilde{\omega}_\psi(v, J_s v) \geq K|v_\mathbb{D}|^2 + |v_X|^2 - 2\|G_s\|\|v_\mathbb{D}\||v_X|,$$

we get tameness as soon as $K$ is bigger than the supremum of $\|G_s\|^2$ on the compact subset $[0,1] \times \mathbb{D} \times U_L$. This proves Lemma 3.4. □

Note also that, because $\pi_\mathbb{R}(U_L) \subset [1, \infty)$, the submanifold $\mathbb{S}^1 \times L \subset \mathbb{D}^2 \times (S\xi \times \mathbb{C})$ is Lagrangian for every $\tilde{\omega}_\psi$. Stokes' formula then ensures that

$$\int_\mathbb{D} \tilde{u}^* \tilde{\omega}_\psi = K\pi$$

for every $u$ in $\mathcal{B}$. If $K$ is sufficiently large to get tameness from Lemma 3.4 and $\tilde{u}$ is $J_s$-holomorphic then $\tilde{u}^* \tilde{\omega}_\psi$ is non-negative on $\mathbb{D}$. Then for any $V \subset \mathbb{D}$, we have

$$\int_V u^*(\omega_\psi \oplus \omega_\mathbb{C}) = \int_V \tilde{u}^* \tilde{\omega}_\psi - \int_V K\omega_\mathbb{D}$$

with both right-hand side integrals in $[0, K\pi]$ so we can choose $A = K\pi$ to finish the proof of Proposition 3.3. □

After those preliminaries, we now consider any sequence $(s_k, u_k)_k$ which has no convergent subsequence in $\mathcal{M}(G)$. The sequence $r_k := \max_\mathbb{D} |du_k|$ is unbounded since otherwise the Arzelà–Ascoli theorem and elliptic regularity would provide a convergent subsequence for $u_k$. Let $z_k$ be a sequence in $\mathbb{D}$ such that $r_k = |du_k(z_k)|$. After passing to a subsequence, we can assume that $z_k$ converges to some $z_\infty$ in $\bar{\mathbb{D}}$ and $r_k$ goes to $+\infty$. We set $\delta_k = d(z_k, \partial\mathbb{D}) = 1 - |z_k|$.

3.1. **Sphere and plane bubbling.** Assume for contradiction that $r_k \delta_k$ is unbounded. (this happens for instance if $z_\infty$ lies in the interior of $\mathbb{D}$). After passing to a subsequence, we can assume that $r_k \delta_k$ is increasing and goes to infinity. We denote by $\mathbb{D}_k$ the open disk with radius $r_k \delta_k$ in $\mathbb{C}$ and consider the map $\Phi_k \colon z \mapsto z_k + z/r_k$ which, due to our choice of $\delta_k$, sends $\mathbb{D}_k$ into $\mathbb{D}$. We set $t_k = \pi_\mathbb{R}(u_k(z_k))$ and

$$v_k \colon \mathbb{D}_k \to S\xi \times \mathbb{C}, \ z \mapsto \tau_{-t_k} \circ u_k \circ \Phi_k.$$

By construction, we have $\sup |dv_k| \leq |dv_k(0)| = 1$ and $\pi_\mathbb{R}(v_k(0)) = 0$. The Arzelà–Ascoli theorem then proves that $v_k$ converges uniformly on compact subsets to some $v \colon \mathbb{C} \to S\xi \times \mathbb{C}$. Since $J$ is translation invariant, we get

$$\bar{\partial}_J v_k = \frac{1}{r_k} d\tau_{-t_k} \circ \bar{\partial}_J u_k(\Phi_k(z))$$

$$= \frac{1}{r_k} d\tau_{-t_k} \circ G_{s_k}(\Phi_k(z), u_k(\Phi_k(z))) \to 0$$

where convergence is uniform on compact sets hence, by elliptic regularity, $v$ is genuinely $J$-holomorphic, and in particular the component $\pi_\mathbb{C} \circ v$ is a classical holomorphic function from $\mathbb{C}$ to $\mathbb{C}$, and $v_\xi := \pi_\xi \circ v$ is a $J_\alpha$-holomorphic map. From Lemma 3.2 it follows that $\pi_\mathbb{C} \circ v$ is bounded, so that this component is in fact constant. In particular we get that, for any compact subset $K$ in $\mathbb{C}$, we have

$$\lim_{k \to \infty} \int_{\Phi_k(K)} u_k^* \omega_\mathbb{C} = 0. \tag{3.2}$$

The component $v_\xi$ by contrast cannot be constant since $|dv(0)| = 1$.

**Lemma 3.5.** *After passing to a subsequence of $u_k$, the sequence $t_k$ diverges towards $-\infty$.*

*Proof.* The $C^0$-estimate from Lemma 3.2 already implies that $t_k = u_k(z_k)$ is bounded above by $T_{\max} = \sup \pi_\mathbb{R}(S\xi \setminus U_+)$ where $U_+$ is given by the lemma. Assume for contradiction that the sequence $t_k$ is bounded below by some $T_{\min}$. In particular $\pi_\mathbb{R} \circ v_\xi$ is bounded above by $T_{\max} - T_{\min}$.

We observe that $v_\xi$ then has finite $\omega$-area. Indeed, for any compact subset $K$ in $\mathbb{C}$ we have

$$\int_K v^* \omega = \lim_{k \to \infty} \int_K v_k^* \omega = \lim_{k \to \infty} \int_{\Phi_k(K)} u_k^* \tau_{-t_k}^* \omega$$

$$= \lim_{k \to \infty} e^{-t_k} \left( \int_{\Phi_k(K)} u_k^* (\omega \oplus \omega_\mathbb{C}) - \int_{\Phi_k(K)} u_k^* \omega_\mathbb{C} \right) \leq A \, e^{-T_{\min}}$$

where the constant $A$ comes from Proposition 3.3 and we used Equation (3.2). Proposition A.1 from the appendix guarantees that either $v_\xi$ is proper or it extends to a $J_\alpha$-holomorphic map from $\mathbb{CP}^1$ to $S\xi$. Since $J_\alpha$ is tamed by the exact symplectic form of $S\xi$ the later possibility is ruled out and $v_\xi$ is proper. Hence the only end of $\mathbb{C}$ is mapped by $v_\xi$ to only one of the two ends of $S\xi$. Since $\pi_\mathbb{R}$ is $J_\alpha$-plurisubharmonic on $S\xi$, the maximum principle implies that $v_\xi$ sends infinity to $S_{+\infty}\xi$, contradicting the above bound by $T_{\max} - T_{\min}$. □

Remember that the Hofer energy of a $J_\alpha$-holomorphic map $v\colon \Sigma \to S\xi$ in the symplectization is defined as
$$E_\alpha(v) = \sup_{\varphi \in \mathcal{F}'} \int_\Sigma v^*\omega_\varphi,$$
where $\mathcal{F}'$ is the space of increasing diffeomorphisms from $\mathbb{R}$ to $(-1,0)$ and $\omega_\varphi = d(e^\varphi \alpha)$ on $S\xi$.

**Lemma 3.6.** *The holomorphic plane $v_\xi$ has finite Hofer energy.*

*Proof.* For any function $\varphi$ in $\mathcal{F}'$ and any $k$, we can choose a function $\psi_k$ in $\mathcal{F}$ that coincides with $\varphi \circ \tau_{-t_k}$ on $(-\infty, 0]$. We know from Lemma 3.5 that $t_k$ diverges towards $-\infty$. Hence, for any radius $r$ there is some $k_0(r)$ such that $t_k \leq -r$ for every $k \geq k_0(r)$ (and $r/r_k < \delta_k$). Since $|du_k| \leq r_k$ and, by definition, $t_k = \pi_\mathbb{R}(u_k(z_k))$ we get by the mean value theorem that $\pi_\mathbb{R}(u_k(z)) \leq t_k + r_k|z - z_k|$ for all $z$. Because $\Phi_k(\mathbb{D}_r) = D(z_k, r/r_k)$, we then get that $\pi_\mathbb{R}(u_k(\Phi_k(\mathbb{D}_r)))$ is contained in $(-\infty, 0]$ for $k \geq k_0(r)$. So our choice of $\psi_k$ gives:
$$\int_{\mathbb{D}_r} v_k^* \omega_\varphi = \int_{\Phi_k(\mathbb{D}_r)} u_k^* \omega_{\varphi \circ \tau_k} = \int_{\Phi_k(\mathbb{D}_r)} u_k^* \omega_{\psi_k}$$
and the later integral is bounded, for $k$ large enough, by $A+1$ where $A$ comes from Proposition 3.3 and where we have used Equation (3.2). We now let $k$ go to infinity then $r$ goes to infinity and finally take the supremum over $\varphi$ to get $E_\alpha(v) \leq A + 1$. □

Since $v_\xi$ has finite Hofer energy, [Hof93, Theorem 31] gives a contractible $T$-periodic Reeb orbit $\gamma$ for $\alpha$ and a sequence $x_k$ such that $v_\xi(e^{2\pi(x_k + iy)})$ converges uniformly to $\gamma(Ty)$ (we cannot hope for convergence without condition on $x_k$ because we haven't made any non-degeneracy assumption on $\alpha$). This contradicts our assumption that $\alpha$ has no contractible closed Reeb orbit so we have proved that $r_k \delta_k$ is bounded.

3.2. **Disk bubbling.** Because $r_k \delta_k$ is bounded, we learn in particular that $z_\infty$ is in $\partial \mathbb{D}$. For notational convenience only, we assume that $z_\infty = 1$. After passing to a subsequence we can assume that $r_k \geq 1$ and $r_k \delta_k$ converges to some non-negative number $\nu$. We set
$$w_k = \left(1 - \frac{1}{r_k}\right)\frac{z_k}{|z_k|}$$
(this extra sequence of points is a minor nuisance needed because when $\nu$ is zero, the naive rescaling could lead to a constant map). We use the rescaling maps
$$\Phi_k(z) = \frac{z + w_k}{1 + \bar{w}_k z}$$
which are automorphisms of $\mathbb{D}$ sending $0$ to $w_k$ and which converge uniformly to the constant map $z \mapsto 1$ on any compact subset of $\mathbb{D}' := \mathbb{D} \setminus \{-1\}$ and have derivative
$$d\Phi_k(z) = \frac{1 - |w_k|^2}{(1 + \bar{w}_k z)^2} dz.$$
In particular:
$$\|d\Phi_k(z)\| = \frac{1 - |w_k|^2}{|1 + \bar{w}_k z|^2} = \frac{2 - 1/r_k}{r_k \mathrm{dist}(-1, \bar{w}_k z)^2}$$
so, for every compact $K \subset \mathbb{D}'$, there are positive constants $C_1(K)$ and $C_2(K)$ such that, for every $z$ in $K$:
$$\frac{C_1(K)}{r_k} \leq \|d\Phi_k(z)\| \leq \frac{C_2(K)}{r_k}.$$
We set $\zeta_k = \Phi_k^{-1}(z_k) = (1 - \delta_k r_k)/(1 + \delta_k r_k - \delta_k) \cdot z_k$ which converges to $\zeta_\infty := (1-\nu)/(1+\nu) \in \mathbb{D}'$. Let $K_0 \subset \mathbb{D}'$ be a compact subset containing all $\zeta_k$ for $k$ large enough.

Our rescaled disk is then $v_k := u_k \circ \Phi_k$ which satisfies: $\|dv_k(z)\| = \|du_k(\Phi_k(z))\| \cdot \|d\Phi_k(z)\|$ for every $z$ since $d\Phi_k(z)$ is an invertible conformal linear map. So for every compact $K \subset \mathbb{D}'$, $\|dv_k\| \leq C_2(K)$ on $K$. In addition each $v_k(\partial \mathbb{D}')$ is in $L$ so $v_k(K)$ is in the $C_2(K)$-neighborhood of $L$ for each convex compact subset $K$. Using an exhaustion of $\mathbb{D}'$ by such subsets we get

that $v_k$ is uniformly bounded on compact subsets of $\mathbb{D}'$. So the Arzelà-Ascoli theorem gives convergence of $v_k$ to some $v\colon (\mathbb{D}', \partial \mathbb{D}') \to (S\xi \times \mathbb{C}, L)$ uniformly on compact subsets of $\mathbb{D}'$. Since $\|dv_k(\zeta_k)\| \geq C_1(K_0)$, we get $dv(\zeta_\infty) \neq 0$ and $v$ is non-constant. Using that $\Phi_k$ is holomorphic we get:
$$\bar{\partial} v_k(z) = G_s(\Phi_k(z), v_k(z)) \circ d\Phi_k(z)$$
which converges to zero uniformly on compact subsets of $\mathbb{D}'$ so, by elliptic regularity, $v$ is $J$-holomorphic.

The energy of $v$ is bounded by Proposition 3.3 since:
$$\int_K v_k^*(\omega \oplus \omega_\mathbb{C}) = \int_{\Phi_k(K)} u_k^*(\omega \oplus \omega_\mathbb{C}) \leq A.$$

Applying the removal of singularity result from the appendix, we can compactify $v$ to a non-constant $J$-holomorphic disk with boundary on $L$. Thus as we wanted to show, $L$ is not weakly exact.

## 4. Proof of the main theorem

**Observation 4.1.** *Additionally to the properties that were proved in* [MNW13] *and recalled in the introduction, the manifold* $(M_\mathbf{k}, \alpha_+)$ *contains a closed pre-Lagrangian submanifold* $P_0$ *such that the restriction of* $\alpha_+$ *to* $P_0$ *is closed and* $\pi_1(P_0)$ *injects into* $\pi_1(M_\mathbf{k})$.

*Proof.* To prove this observation, let us briefly sketch how the manifolds $(M_\mathbf{k}, \alpha_\pm)$ were constructed. Define $\overline{M}$ to be $\mathbb{R}^l \times \mathbb{R}^{l+1}$ with coordinates $(t_1, \ldots, t_l; \theta_0, \ldots, \theta_l)$, and choose on $\overline{M}$ the two contact forms
$$\alpha_\pm := \pm e^{t_1 + \cdots + t_l} d\theta_0 + e^{-t_1} d\theta_1 + \cdots + e^{-t_l} d\theta_l.$$
The desired manifold $M_\mathbf{k}$ is a quotient of $\overline{M}$ by a group of transformations depending on $\mathbf{k}$.

There are two natural group actions on $\overline{M}$ that preserve the contact forms $\alpha_\pm$. The most obvious one is the one of $\mathbb{R}^{l+1}$ by translations on the second factor, i.e.
$$(t_1, \ldots, t_l; \theta_0, \ldots, \theta_l) + (\varphi_0, \ldots, \varphi_l) := (t_1, \ldots, t_l; \theta_0 + \varphi_0, \ldots, \theta_l + \varphi_l).$$
If we choose any lattice $\Lambda' \subset \mathbb{R}^{l+1}$, the quotient $\overline{M}/\Lambda'$ will be diffeomorphic to $\mathbb{R}^l \times \mathbb{T}^{l+1}$, and both $\alpha_+$ and $\alpha_-$ descend to this quotient.

We can also define an action of $\mathbb{R}^l$ on $\overline{M}$ by
$$(\tau_1, \ldots, \tau_l) \cdot (t_1, \ldots, t_l; \theta_0, \theta_1, \ldots, \theta_l) :=$$
$$(t_1 + \tau_1, \ldots, t_l + \tau_l; e^{-(\tau_1 + \cdots + \tau_l)} \theta_0, e^{\tau_1} \theta_1, \ldots, e^{\tau_l} \theta_l).$$
This second action also preserves the contact forms $\alpha_\pm$, but it does not commute with the first action, hence it is not obvious how to combine it with the translations by $\Lambda'$ to produce a compact quotient of $\overline{M}/\Lambda'$. In other word, it not immediate to find a lattice in the relevant semi-direct product $\mathbb{R}^l \ltimes \mathbb{R}^{l+1}$. Nonetheless according to [MNW13, Lemma 8.3], there exist lattices $\Lambda \subset \mathbb{R}^l$, and $\Lambda' \subset \mathbb{R}^{l+1}$ such that the restriction of the $\mathbb{R}^l$-action to $\Lambda$ preserves $\mathbb{R}^l \times \Lambda'$. They depend on $\mathbf{k}$ but this dependence will be suppressed for easier reading. It follows that there is a well-defined action of $\Lambda$ on $\overline{M}/\Lambda'$, so that $M := (\overline{M}/\Lambda')/\Lambda$ is a smooth manifold, and since the $\alpha_\pm$ are invariant under both the $\Lambda$- and the $\Lambda'$-action, they induce contact forms on $M$. The projection $\mathbb{R}^l \times \mathbb{R}^{l+1} \to \mathbb{R}^l$ induces on $M$ the structure of a $\mathbb{T}^{l+1}$-bundle over $\mathbb{T}^l$.

Consider any of the fibers $\{\mathbf{t}\} \times (\mathbb{R}^{l+1}/\Lambda')$ in $\overline{M}/\Lambda'$ with $\mathbf{t} \in \mathbb{R}^l$. This fiber is a torus which is pre-Lagrangian in $(\overline{M}/\Lambda', \alpha_+)$, because $\alpha_+$ restricts to a constant 1-form on it. Clearly this torus embeds in $M$ under the projection $\overline{M}/\Lambda' \to M$. We choose for $P_0$ the image of this embedding. By construction, $\pi_1(P_0) = \Lambda'$ embeds into $\pi_1(M) = \Lambda \ltimes \Lambda'$. □

*Proof of Theorem 1.1.* We set $V_\mathbf{k} = \mathbb{T}^2 \times M_\mathbf{k}$. The pre-Lagrangian submanifold $P_0 \subset (M_\mathbf{k}, \alpha_+)$ from Observation 4.1 extends to a pre-Lagrangian submanifold $P := \{0\} \times \mathbb{S}^1 \times P_0$ in $(V_\mathbf{k}, \xi_n)$, because the restriction of
$$\lambda_n = \frac{1 + \cos(ns)}{2} \alpha_+ + \frac{1 - \cos(ns)}{2} \alpha_- + \sin(ns)\, dt$$

to $\{0\} \times \mathbb{S}^1 \times P_0$ is the closed 1-form $\alpha_+|_{TP_0}$.

The contactomorphism $\Psi_{n,m} \colon (s,t,\theta) \mapsto (s + 2\pi m/n, t, \theta)$ obviously displaces $P$ from itself so we only need to check that $P$ is weakly exact and apply Theorem 2.3 to get that $\Psi_{n,m}$ is not symplectically pseudo-isotopic to the identity.

We will now show $\pi_2(V_{\mathbf{k}}, P) = 0$ which implies that $P$ is weakly exact. Because $\pi_1(P_0)$ embeds into $\pi_1(M_{\mathbf{k}})$ we get that $\pi_1(P)$ embeds into $\pi_1(V_{\mathbf{k}})$. The long exact sequence of the pair $(V_{\mathbf{k}}, P)$ contains

$$\pi_2(V_{\mathbf{k}}) \to \pi_2(V_{\mathbf{k}}, P) \to \pi_1(P) \to \pi_1(V_{\mathbf{k}})$$

where $\pi_2(V_{\mathbf{k}}) = 0$ (because by construction the universal cover of $M_{\mathbf{k}}$ is a Euclidean space) and $\pi_1(P) \hookrightarrow \pi_1(V_{\mathbf{k}})$ so $\pi_2(V_{\mathbf{k}}, P) = 0$.

We now prove the last part of the main theorem, about conjugations. Consider the family of contactomorphisms

$$\Phi_\tau \colon (\mathbb{T}^2 \times M_{\mathbf{k}}, \xi_n) \to (\mathbb{T}^2 \times M_{\mathbf{k}}, \xi_n),\ (s,t;x) \mapsto (s, t + \tfrac{2\pi m \tau}{n}; x)$$

for $\tau \in [0,1]$.

The contactomorphism $\Phi_1$ is conjugated to $\Psi_{n,m}$ by the diffeomorphism

$$A \colon (V_{\mathbf{k}}, \xi_n) \to (V_{\mathbf{k}}, \xi_n),\ (s,t;x) \mapsto (t, -s; x)$$

which satisfies $A \Phi_1 A^{-1} = \Psi_{n,m}$. But $A$ cannot be replaced by a contactomorphism since $\Phi_1$ is, by construction, contact isotopic to the identity. $\square$

## Appendix A. Removal of singularities

Even though removal of singularities for holomorphic curves is a recurring topic, we haven't been able to find a reference proving precisely the following versions needed here (see e.g. [MS04, Section 4.1]). We use the notations $\mathbb{D} = \{z \in \mathbb{C}\mid |z| \leq 1\}$ and $\mathbb{D}^+ = \{z \in \mathbb{D}\mid \operatorname{Im}(z) \geq 0\}$

**Proposition A.1.** *Let $(W, \omega)$ be a symplectic manifold that does not need to be compact or geometrically bounded in any sense, and let $J$ be a compatible almost complex structure on $W$.*

(a) *Assume that $L \subset W$ is a properly embedded (but not necessarily closed) Lagrangian submanifold. Every $J$-holomorphic map*

$$(\mathbb{D}^+ \setminus \{0\}, [-1,1] \setminus \{0\}) \to (W, L)$$

*whose $\omega$-area is finite is either proper or extends to a $J$-holomorphic map from $(\mathbb{D}^+, [-1,1])$ to $(W, L)$. In particular if $L$ is compact then $u$ extends.*

(b) *Every $J$-holomorphic map*

$$\mathbb{D} \setminus \{0\} \to W$$

*whose $\omega$-area is finite is either proper or extends to a $J$-holomorphic map from $\mathbb{D}$ to $W$.*

All distances are measured using the Riemannian metric $g := \omega(\cdot, J\cdot)$, and the $\omega$-area of any $u \colon \Sigma \to W$ is given by $E(u) = \int_\Sigma u^*\omega$.

To treat case (a) and (b) in a unified form, we introduce the following notations and conventions. Let $\Omega$ denote either $\mathbb{D}$ or $\mathbb{D}^+$ and $\Omega^* = \Omega \setminus \{0\}$. Let $\partial_L \Omega$ denote $[-1,1]$ if $\Omega = \mathbb{D}^+$ and the empty set otherwise. We also set $\partial_L \Omega^* = \partial_L \Omega \setminus \{0\}$. In case (b) it is understood that $L$ is empty. With those conventions, $u$ is always a map of pairs from $(\Omega^*, \partial_L \Omega^*)$ to $(W, L)$. Finally, denote the circle or half-circle $\{z \in \Omega \mid |z| = r\}$ by $\Gamma_r$.

The following observation which will be important for our argument is a simple topological remark that does not depend on $J$-holomorphicity.

**Lemma A.2.** *Let $u$ be a continuous map $\Omega^* \to W$ that is not proper and does not extend continuously over $0$. Then there are two sequences $(z_k)_k$ and $(w_k)_k$ in $\Omega^*$ both converging to $0$ such that $u(z_k)$ and $u(w_k)$ both converge but $\lim u(z_k) \neq \lim u(w_k)$.*

*Proof.* Since $u$ is not proper there exists a compact subset $K \subset W$ such that $\Omega_K := u^{-1}(K)$ is *not* compact (and not empty). This implies that $\Omega_K$ contains a sequence $(z_k)_k$ that has no converging subsequence, but since $\Omega_K$ is closed in $\Omega^*$, it follows that $z_k \to 0$ and then, after passing to a subsequence, also that $u(z_k)$ converges to some $p \in K$.

Since by assumption $u$ does not extend continuously over $0$, there has to be a second sequence $(z'_k)_k$ in $\Omega^*$ also converging to $0$ such that its image $u(z'_k)$ does *not* converge to $p$, even after passing to a subsequence. We can take a small ball $B_\varepsilon(p) \subset W$ of radius $\varepsilon > 0$ which is compact and such that, after passing to a subsequence, all $u(z'_k)$ are outside $B_\varepsilon(p)$. There is a sequence of continuous arcs $\gamma_k \colon [0,1] \to \Omega^*$ from $z_k$ to $z'_k$ that converge uniformly to $0$. By continuity the image of $u \circ \gamma_k$ needs to intersect $\partial B_\varepsilon(p)$, and hence gives rise to a sequence $(w_k)_k \subset \Omega^*$ with $w_k \to 0$ and $u(w_k) \in \partial B_\varepsilon(p)$. The later is compact hence a subsequence of $u(w_k)$ converges to some $q \in \partial B_\varepsilon(p)$. $\square$

The next lemma crucially uses the monotonicity lemma on compact subsets.

**Lemma A.3.** *Let $u \colon (\Omega^*, \partial_L \Omega^*) \to (W, L)$ be a $J$-holomorphic curve with finite $\omega$-area. Assume there is a sequence $(z_k)_k$ in $\Omega^*$ with $z_k \to 0$ and $u(z_k) \to p \in W$. For every positive distance $\delta$ there is a radius $r_\delta$ such that every arc $u(\Gamma_r)$ with $r \in (0, r_\delta]$ intersects $B_\delta(p)$.*

*Proof.* Denote the ball $B_\delta(p)$ by $K$. Shrinking $\delta$ if necessary, we may assume that $K$ (and its intersection with $L$) are compact. According to the monotonicity lemma (see e.g. [Sik94, Propositions 4.3.1 and 4.7.2]) there are positive constants $r_K$ and $C_K$ such that for every $r < r_K$, and for every $x \in K$, and for every non-constant holomorphic map
$$v \colon (G, \partial G) \to \Big(B_r(x), \partial B_r(x) \cup (L \cap B_r(x))\Big)$$
defined on a compact domain $G \subset \mathbb{C}$ with piecewise smooth boundary, and such that $x \in v(G)$, we have
$$E(v) \geq C_K r^2$$
(recall that $L$ is empty in case (b)).

Assume for contradiction that no $r_\delta$ is small enough to get the announced conclusion. It means there is a sequence $r_j \to 0$ such that none of the $u(\Gamma_{r_j})$ intersects $B_\delta(p)$. Let $\delta'$ be $\min\{r_K, \delta/2\}$. After passing to subsequences from both $z_k$ and $r_j$ we can assume $r_{k+1} < |z_k| < r_k$ and $p_k := u(z_k) \in B_{\delta'}(p)$ for every $k$.

Let now $U_k \subset \Omega^*$ be the connected component of $u^{-1}(B_\delta(p))$ containing the point $z_k$. Since the arcs $\Gamma_{r_k}$ disconnect $\Omega^*$, all the $U_k$'s are pairwise disjoint. Using Sard's theorem we find for every $k$, a radius $\delta_k \in [\frac{\delta'}{2}, \delta']$ such that $u^{-1}(B_{\delta_k}(p_k)) \cap U_k$ has piecewise smooth boundary. This allows us to apply the monotonicity lemma quoted above to each $U_k$ to see that $\int_{U_k} u^*\omega$ is at least $C_K(\delta'/2)^2$, but this contradicts the finiteness of $\int_{\Omega^*} u^*\omega$. $\square$

The following lemma is classical and can be found for example in the proof of [Sik94, Theorem 4.5.1].

**Lemma A.4.** *Let $u \colon (\Omega^*, \partial_L \Omega^*) \to (W, L)$ be a $J$-holomorphic curve with finite $\omega$-area. For any $\varepsilon > 0$, there are arbitrary short arcs $u(\Gamma_r)$ with $r \in (0, \epsilon]$.*

*Proof.* Suppose the statement were wrong, then there would exist a positive $\varepsilon$ such that every arc $u(\Gamma_r)$ with $r \in (0, \varepsilon]$ will be longer than some constant $\delta > 0$. We have $\omega(\partial_r u, \partial_\theta u) = -\omega(\partial_\theta u, \partial_r u) = \omega(\partial_\theta u, \frac{1}{r} J \cdot \partial_\theta u) = \frac{1}{r} g(\partial_\theta u, \partial_\theta u) = \frac{1}{r}|\partial_\theta u|^2$. Setting $N = 1$ in case (a) and $N = 2$ in case (b), we compute
$$E(u) = \int_{\Omega^*} \frac{1}{r} \Big|\frac{\partial u}{\partial \theta}\Big|^2 dr \wedge d\theta \geq \int_0^\varepsilon \Big(\frac{1}{r} \int_0^{N\pi} \Big|\frac{\partial u}{\partial \theta}\Big|^2 d\theta\Big) dr = \int_0^\varepsilon \frac{1}{r} \Big\|\frac{\partial u}{\partial \theta}(r, \cdot)\Big\|_{L^2([0, N\pi])}^2 dr$$
$$= \frac{1}{N\pi} \int_0^\varepsilon \frac{1}{r} \|1\|_{L^2}^2 \Big\|\frac{\partial u}{\partial \theta}(r, \cdot)\Big\|_{L^2}^2 dr \geq \int_0^\varepsilon \frac{1}{N\pi r} \Big\langle 1, \Big|\frac{\partial u}{\partial \theta}\Big|\Big\rangle_{L^2}^2 dr$$
$$= \int_0^\varepsilon \frac{1}{N\pi r} \Big\|\frac{\partial u}{\partial \theta}\Big\|_{L^1}^2 dr = \int_0^\varepsilon \frac{1}{N\pi r} \ell(u(\Gamma_r))^2 dr \geq \int_0^\varepsilon \frac{\delta^2}{N\pi r} dr = \infty \ .$$

This contradicts our assumption that $E(u)$ is finite. $\square$

*Proof of Proposition A.1.* Let $u\colon (\Omega^*, \partial_L\Omega^*) \to (W, L)$ be a $J$-holomorphic curve that has finite $\omega$-area and is not proper. Assume for contradiction that $u$ does not extend continuously over the origin.

Lemma A.2 gives sequences $z_k$ and $w_k$ in $\Omega^*$ converging to 0 such that $u(z_k)$ converges to some $p$ and $u(w_k)$ converges to some $q$ different from $p$. Lemma A.3 applied to these sequences with $\delta = d(p,q)/3$ gives a radius $r_\delta$ such that every arc $u(\Gamma_r)$ with $r \leq r_\delta$ intersects both $B_\delta(p)$ and $B_\delta(q)$. Since the distance $d(B_\delta(p), B_\delta(q))$ is at least $\delta$, the length of each arc $u(\Gamma_r)$ needs to be larger than $\delta > 0$, but this is a contradiction to Lemma A.4. Hence $u$ admits a continuous extension over the origin. In particular $u$ has a relatively compact image thus we can apply the usual discussion of removal of singularities to prove that the extension is $J$-holomorphic. This can be done either by proving that $u$ is $W^{1,p}$ for some $p > 2$ and use elliptic regularity bootstrapping or by proving that $du$ can also be seen as a finite area pseudo-holomorphic curve, hence enabling a geometric bootstrapping argument, see the extensive discussion in [FQ15]. For the purposes of the present paper, none of this extra work is required, see Remark A.5. $\square$

*Remark* A.5. The arguments in this appendix only show that if a $J$-holomorphic plane $u\colon \mathbb{C} \to W$ with finite positive area $A$ is not proper then it needs to extend to a continuous map $\hat{u}\colon \hat{\mathbb{C}} \to W$. To obtain smoothness of $\hat{u}$, we had to invoke elliptic regularity which even told us that $\hat{u}$ is $J$-holomorphic.

In fact in the setup of this article, it suffices to show that $(W, \omega)$ contains a smooth sphere with positive $\omega$-area (so that in particular $\omega$ is not exact). This weaker statement can be easily obtained directly with the proofs given in this appendix. The key point is Lemma A.4 which proves the existence of an increasing sequence of radii $R_k$ diverging to $+\infty$ such that the lengths $l_k$ of the circles $u(\Gamma_{R_k})$ converge to zero. If $k$ is large enough then those circles will lie in some ball around $\hat{u}(\infty)$. We can then fill them by smooth disks whose Euclidean area is controlled by $c_0 l_k^2$ (see for instance [Hum97, Appendix A]) and hence converge to zero. Because $\omega$ is continuous, its integral over those disks will also converge to zero. So, for any positive $\varepsilon$, there is a piecewise smooth $v$ from $\hat{\mathbb{C}}$ to $W$ which coincides with $u$ on a very large disk and whose $\omega$-area is at least $E(u) - \varepsilon$. Then one can smooth $v$ without changing its $\omega$-area (here Stokes' theorem for manifolds with piecewise smooth boundary is enough). The same strategy applies for disks in order to prove that a Lagrangian submanifold is not weakly exact.

## References

[ALP94]  M. Audin, F. Lalonde, and L. Polterovich, *Symplectic rigidity: Lagrangian submanifolds*, Holomorphic curves in symplectic geometry, Progr. Math., vol. 117, Birkhäuser, Basel, 1994, pp. 271–321.

[CE12]  K. Cieliebak and Y. Eliashberg, *From Stein to Weinstein and back*, American Mathematical Society Colloquium Publications, vol. 59, American Mathematical Society, Providence, RI, 2012, Symplectic geometry of affine complex manifolds.

[EHS95]  Y. Eliashberg, H. Hofer, and D. Salamon, *Lagrangian intersections in contact geometry*, Geom. Funct. Anal. **5** (1995), no. 2, 244–269.

[FQ15]  U. Fuchs and L. Qin, *A theorem on the removal of boundary singularities of pseudo-holomorphic curves*, J. Symplectic Geom. **13** (2015), no. 3, 527–544.

[GM15]  E. Giroux and P. Massot, *On the contact mapping class group of Legendrian circle bundles*, arXiv:1506.01162, 2015.

[Gom98]  R. Gompf, *Handlebody construction of Stein surfaces*, Ann. of Math. (2) **148** (1998), no. 2, 619–693.

[Gro85]  M. Gromov, *Pseudo holomorphic curves in symplectic manifolds*, Invent. Math. **82** (1985), 307–347.

[Hof93]  H. Hofer, *Pseudoholomorphic curves in symplectizations with applications to the Weinstein conjecture in dimension three*, Invent. Math. **114** (1993), no. 3, 515–563.

[Hum97]  C. Hummel, *Gromov's compactness theorem for pseudo-holomorphic curves*, Progress in Mathematics, vol. 151, Birkhäuser Verlag, Basel, 1997.

[MNW13]  P. Massot, K. Niederkrüger, and C. Wendl, *Weak and strong fillability of higher dimensional contact manifolds*, Invent. Math. **192** (2013), no. 2, 287–373.

[MS04]  D. McDuff and D. Salamon, *J-holomorphic curves and symplectic topology*, Colloquium Publications. American Mathematical Society 52. Providence, RI: American Mathematical Society (AMS)., 2004.

[Sik94]  J.-C. Sikorav, *Some properties of holomorphic curves in almost complex manifolds*, Holomorphic curves in symplectic geometry, Progr. Math., vol. 117, Birkhäuser, Basel, 1994, pp. 165–189.

(P. Massot) Centre de Mathématiques Laurent Schwartz, École Polytechnique, 91128 Palaiseau Cedex, FRANCE

*E-mail address*, P. Massot: `patrick.massot@polytechnique.edu`

(K. Niederkrüger) Alfréd Rényi Institute of Mathematics, Hungarian Academy of Sciences, POB 127, H-1364 Budapest, HUNGARY

(K. Niederkrüger) Institut de mathématiques de Toulouse, Université Paul Sabatier – Toulouse III, 118 route de Narbonne, F-31062 Toulouse Cedex 9, FRANCE

*E-mail address*: `niederkr@math.univ-toulouse.fr`